\documentclass{article}
\usepackage{amsmath, amssymb}
\usepackage{amsthm}
\usepackage{enumerate}
\newtheorem{theorem}{Theorem}
\newtheorem{lemma}{Lemma}

\def\Todorcevic{Todor\v{c}evi\'{c}}
\def\Lucke{L\"{u}cke}
\DeclareMathOperator{\hht}{ht}
\DeclareMathOperator{\cf}{cf}
\DeclareMathOperator{\otp}{otp}
\DeclareMathOperator{\dom}{dom}

\DeclareMathOperator{\Lim}{Lim}

\makeatletter
\newcommand{\fcolon}{%
  \mathrel{\mathpalette\fcolon@\relax}%
}
\newcommand{\fcolon@}[2]{%
  \sbox\z@{$\m@th#1:$}%
  \vbox to\ht\z@{%
    \hbox{$\m@th#1.$}%
    \vss
    \hbox{$\m@th#1.$}%
    \vss
    \hbox{$\m@th#1.$}%
  }%
}
\makeatother

\title{An Inconsistent Forcing Axiom at \(\omega_2\)}
\author{Stevo \Todorcevic \and Shihao Xiong}

\newcommand{\Addresses}{{
  \bigskip
  \footnotesize

  Stevo \Todorcevic, \textsc{Department of Mathematics, University of Toronto, Canada},
  \textit{E}: \texttt{stevo@math.toronto.edu}\par\nopagebreak \textsc{Institut de Math\'ematiques de Jussieu, Paris, France},
  \textit{E}: \texttt{stevo.todorcevic@imj-prg.fr}

  \medskip

  Shihao Xiong, \textsc{Department of Mathematics, Cornell University, Ithaca, NY, United States},
  \textit{E}: \texttt{sx76@cornell.edu}
}}

\date{August 3, 2020}

\begin{document}
\maketitle
\begin{abstract}
We show that the forcing axiom for countably compact, \(\omega_2\)-Knaster, well-met posets is inconsistent. This is supplemental to an inconsistency result of Shelah~\cite{shelah1982generalized} and sets a new limit to the generalization of Martin's Axiom to the stage of \(\omega_2\).
\end{abstract}

\section{Introduction}
As an attempt to generalize Martin's Axiom to the stage of \(\omega_2\), Baumgartner~\cite{baumgartner1983iterated} proved the following theorem which shows the consistency of a forcing axiom for \(\omega_1\)-linked posets.

\begin{theorem}[Baumgartner~\cite{baumgartner1983iterated}]\label{BA}
It is consistent with CH and \(2^{\omega_1}>\omega_2\) that \(FA_\kappa(\mathfrak{B})\) holds for every \(\kappa < 2^{\omega_1}\), in which \(\mathfrak{B}\) is the class of all \(\sigma\)-closed, \(\omega_1\)-linked and well-met posets.
\end{theorem}

Here a poset \(P\) is \textit{well-met} if every pair of compatible conditions from \(P\) has a greatest lower bound. \(FA_\kappa(P)\) is the statement that for any family, \(\{D_\xi : \xi < \kappa\}\), of dense open subsets of \(P\), there exists a sufficiently generic filter \(F\subset P\) which intersects every \(D_\xi\). \(FA_\kappa(\mathfrak{B})\) states that \(FA_\kappa(P)\) holds for every \(P\in \mathfrak{B}\).

Many applications of Baumgartner's forcing axiom was studied by Tall~\cite{tall1994some}. However, Shelah~\cite{shelah1982generalized} showed that the full generalization of Martin's Axiom is inconsistent with CH + \(2^{\omega_1}>\omega_2\).

\begin{theorem}[Shelah~\cite{shelah1982generalized}]\label{ccInconsistent}
If CH holds and \(2^{\omega_1}>\omega_2\), then there is a \(\sigma\)-closed poset \(P\) satisfying the \(\omega_2\)-chain conditions such that FA\(_{\omega_2}(P)\) fails.
\end{theorem}

The inconsistency of such a na\"{i}ve generalization of Martin's Axiom does not come from the chain condition directly. Actually the counterexample constructed by Shelah is \(\omega_1\)-linked and the inconsistency comes from the absence of the well-met condition.

The first author of this paper observed that even if we add the well-metness condition, the forcing axioms can still fail to be consistent due to an excessively strong chain condition. Namely, we show the following result.

A poset \(P\) is \(\kappa\)-\textit{Knaster} if for every \(S\subset P\) with cardinality \(\kappa\), there exists a subset \(S_0 \subset S\) with cardinality \(\kappa\) such that the elements in \(S_0\) are pairwisely compatible in \(P\). A poset is \textit{countably compact} if every countable decreasing chain from the poset has a greatest lower bound. Let \(\mathfrak{K}\) be the class of all countably compact, \(\omega_2\)-Knaster and well-met posets.

\begin{theorem}\label{FAKnaster}
\(FA_{\omega_2}(\mathfrak{K})\) is inconsistent with ZFC + CH + \(2^{\omega_1}>\omega_2\).
\end{theorem}

\section{(S) Implies \(\Box_{\omega_1}\)}
We first introduce Shelah's forcing axiom (S)~\cite{shelah1978weak}.

A poset \(P\) is \(\omega_2\)-\textit{normal} if for any sequence \( \langle p_\xi : \xi\in \omega_2\rangle \) of conditions from \(P\), there exists a club \(C\subset \omega_2\) and a regressive function \(f:\omega_2\rightarrow\omega_2\) such that for any \(\alpha, \beta \in C\) with \( \cf\alpha, \cf\beta>\omega \) and \(f(\alpha) = f(\beta)\), \(p_\alpha\) and \(p_\beta\) are always compatible.

\begin{theorem}[Shelah~\cite{shelah1978weak}]\label{ShelahFA}
It is consistent with CH + \(2^{\omega_1}>\omega_2\) that \(FA_\kappa(\mathfrak{S})\) holds for every \(\kappa < 2^{\omega_1}\), in which \(\mathfrak{S}\) is the class of all \(\omega_2\)-normal, countably compact and well-met posets.
\end{theorem}

Let (S) denote Shelah's forcing axiom \(\forall \kappa < 2^{\omega_1}, FA_\kappa(\mathfrak{S})\). Clearly, (S) is a direct consequence of our inconsistent axiom \(FA_{\omega_2}(\mathfrak{K})\). We show that (S) implies \(\Box_{\omega_1}\).

\begin{lemma}\label{ShFAImpliesSquare}
(S) + CH + \(2^{\omega_1}>\omega_2\) implies \(\Box_{\omega_1}\).
\end{lemma}
\begin{proof}
We study the poset \(P\) of all partial functions \(p \fcolon S^2_1 \rightarrow P(\omega_2)\) in which \(S^2_1 = \{\alpha\in\omega_2: \cf\alpha = \omega_1\}\), satisfying:
\begin{enumerate}
     \item\label{ShFAImpliesSquare-Condition-1} \(\dom p\) is countable.
     \item\label{ShFAImpliesSquare-Condition-2} For \(\alpha\in\dom p\), \(p(\alpha)\) is a countable subset of \(\alpha\) and closed in \(\alpha\).
     \item\label{ShFAImpliesSquare-Condition-3} If \(\delta<\gamma\) are both in \(\dom p\) and \(\xi\) is a limit point of both \(p(\delta)\) and \(p(\gamma)\), then \(p(\delta)\cap \xi = p(\gamma)\cap \xi\).
\end{enumerate}
  
The order of the poset is defined by setting \(p\leq q\) if and only if
\begin{enumerate}
     \setcounter{enumi}{3}
     \item\label{ShFAImpliesSquare-Condition-4} \(\dom q\subset \dom p\).
     \item\label{ShFAImpliesSquare-Condition-5} For \(\alpha\in \dom q\), \(q(\alpha) \sqsubseteq p(\alpha)\).\footnote{\(q(\alpha) \sqsubseteq p(\alpha)\) means that \(p(\alpha)\) is an end extension of \(q(\alpha)\).}
     \item\label{ShFAImpliesSquare-Condition-6} For \(\alpha\in \dom q\),
     \[\min p(\alpha) \setminus q(\alpha) \geq \sup\bigcup_{\xi\in\dom q} (q(\xi)\cap\alpha), \]
     and \[\min p(\alpha) \setminus q(\alpha) \geq \sup(\dom q\cap \alpha).\]
\end{enumerate}

Clearly, this order is a transitive one. We first verify that \(P\) is well-met. Given \(p, q\in P\) of which the compatibility is witnessed by \(r_0\), define \(r : \dom p \cup \dom q \rightarrow P(\omega_2)\) such that \(r(\alpha) = p(\alpha)\cup q(\alpha)\). Clearly every \(r(\alpha)\) is a countable closed subset of \(\alpha\) and is an initial segment of \(r_0(\alpha)\). For any \(\delta<\gamma \in \dom r\) and \(\xi\), a limit point of both \(r(\delta)\) and \(r(\gamma)\), it must be true that \(\xi\) is also a limit point of both \(r_0(\delta)\) and \(r_0(\gamma)\). Since \(r_0\) is a condition from \(P\), we have that \(\xi\cap r(\delta)= \xi\cap r_0(\delta) = \xi\cap r_0(\gamma)= \xi\cap r(\gamma)\). Therefore, \(r\) is a condition from \(P\). Meanwhile, \(r\) extends both \(p\) and \(q\). For any \(\alpha\in \dom r\), since both \(p(\alpha)\) and \(q(\alpha)\) are initial segments of \(r_0(\alpha)\), \(r(\alpha)\) must be an end-extension of both of them. Condition~\ref{ShFAImpliesSquare-Condition-6} is witnessed by \(r_0\) and obviously \(r\) is the greatest lower bound of \(p\) and \(q\).
  
To see that \(P\) is countably compact, let \(\langle p_i: i\in\omega\rangle\) be a countable decreasing chain from \(P\). Define \(r : \bigcup_{i\in\omega} \dom p_i \rightarrow P(\omega_2)\) such that \(r(\alpha) = cl(\bigcup_{i\in\omega} p_i(\alpha))\). It is clear that \(r(\alpha)\) must be a countable subset of \(\alpha\) since there is at most one point in \(r(\alpha) - \bigcup_{i\in\omega} p_i(\alpha)\), which is \(\sup (\bigcup_{i\in\omega}p_i(\alpha))\) if it is not already in some \(p_i(\alpha)\). Let's say that this point is \textit{special} for \(\alpha\).
  
It suffices to verify condition~\ref{ShFAImpliesSquare-Condition-3} for \(r\). Let \(\xi\) be a limit point of both \(r(\gamma)\) and \(r(\delta)\) and \(i\) be the minimal natural number such that both \(\gamma\) and \(\delta\) are in \(\dom p_i\). Without loss of generality, we assume that \(\gamma < \delta\). If \(\xi\) is not special for either \(\gamma\) or \(\delta\), there exists \(j\in\omega\) such that \(\xi\in p_j(\gamma)\cap p_j(\delta)\). Therefore, \(\xi\cap r(\gamma) = \xi\cap p_j(\gamma)= \xi\cap p_j(\delta)= \xi\cap r(\delta)\). If \(\xi\) is special for \(\delta\), there exists \(j>i\) such that we can find an \(\eta\in p_j(\delta)\setminus p_i(\delta)\) and \(\eta<\xi<\gamma\). This contradicts condition~\ref{ShFAImpliesSquare-Condition-6}. If \(\xi\) is special for \(\gamma\) but not special for \(\delta\), there exists \(j_0\geq i\) such that \(\xi\in p_{j_0}(\delta)\). Now for any \(j>j_0\), by condition~\ref{ShFAImpliesSquare-Condition-6}, \(\min p_j(\gamma)\setminus p_{j_0}(\gamma) \geq \xi\) which is a contradiction.
  
Therefore, \(r\) extends the sequence \(\langle p_i: i\in\omega\rangle\) and it must be the greatest lower bound of this sequence by definition.
  
By simple density arguments, we can see that a sufficiently generic object must be a sequence \(C = \langle C_\beta: \beta\in S^2_1 \rangle\) such that every \(C_\beta\) is a club subset of \(\beta\) with ordertype \(\omega_1\).  We can extend the definition to the ordinals with countable cofinality. For \(\alpha\in\omega_2\) such that \(\cf\alpha=\omega\), if there is some \(\beta\in S^2_1\) such that \(\alpha\) is a limit point of \(C_\beta\), then define \(C_\alpha = C_\beta\cap\alpha\). By the coherence of \(C\) given by condition~\ref{ShFAImpliesSquare-Condition-6}, it does not depend on our choice of \(\beta\). If there is no such \(\beta\), then just define \(C_\alpha\) to be any cofinal subset with ordertype \(\omega\). This extension must be a \(\Box_{\omega_1}\)-sequence.
  
Now it suffices to verify that \(P\) is \(\omega_2\)-normal. Fix any sequence \(\langle p_\xi : \xi\in\omega_2 \rangle\) from \(P\) and a club \(C \subset \omega_2\) such that \(\dom p_\xi \subset \eta\) for any \(\xi < \eta \in C\).
  
Given \(\xi\in S^2_1\), define \[s_\xi = \sup ( (\dom p_\xi \cup \bigcup_{\alpha\in\dom p_\xi}p_\xi(\alpha)) \cap \xi) + 1.\] \(s_\xi\) is always strictly less than \(\xi\) since \(\cf \xi = \omega_1\).
  
For every \(\xi\in S^2_1\), let \(\gamma_0 = \otp(\dom p_\xi \cap s_\xi)\) and \(\gamma_1 = \otp(\dom p_\xi \setminus s_\xi)\). Without loss of generality, we may assume that both of them are independent from the choice of \(\xi\). Now we fix a function \(F:\omega_2\times \omega_1 \rightarrow V \) such that \[\forall \iota \in \omega_2, F(\iota, \cdot):\omega_1\rightarrow V \text{ emumerates } (\iota \times [\iota]^\omega)^{\gamma_0} \times ([\iota]^\omega)^{\gamma_1}.\] The existence of such a function \(F\) is given by our assumption CH.
  
For \(\xi\in\omega_2\), define \(f(\xi) = \langle s_\xi, \tau_\xi \rangle\) in which \(\tau_\xi\) is the unique \(\tau\) such that the first term of \(F(s_\xi, \tau)\) codes the restriction of the function \(p_\xi\) on \(\xi\), i.e., \[\{ \langle \alpha, p_\xi(\alpha) \rangle :  \alpha\in\dom p_\xi \cap \xi\},\] and the second term of \(F(s_\xi, \tau)\) codes the function \(\xi_i \mapsto p_\xi(\xi_i)\cap s_\xi\) in which \(i\in \gamma_1\) and \(\xi_i\) is the \(i\)-th element of \(\dom p_\xi \setminus s_\xi\).
  
\(f\) is a desired regressive function such that whenever \(\xi < \eta\) are in \(C\) with cofinality \(\omega_1\) and \(f(\xi) = f(\eta)\), \(p_\xi\) is compatible with \(p_\eta\). Given \(\xi<\eta\) with such properties, we define the condition \(r\) by simply taking \(r(\beta) = p_\xi(\beta) \cup p_\eta(\beta)\) for any \(\beta\in \dom p_\xi \cup \dom p_\eta\). Note that either \(p_\xi(\beta) = p_\eta(\beta)\) or one of them is not defined, in which case we understand it as an empty set. Therefore, \(r\) extends both \(p_\xi\) and \(p_\eta\) as long as it is a condition.
  
It suffices to verify condition~\ref{ShFAImpliesSquare-Condition-3} for \(r\). Let \(s = s_\xi = s_\eta\). Assume that \(\alpha\) is a limit point of both \(r(\gamma)\) and \(r(\delta)\). We may assume that \(\gamma\in \dom p_\xi \setminus s\) and \(\delta\in \dom p_\eta \setminus s\). By the definition of the club \(C\), \(\gamma<\eta\leq \delta\), so we have \(\alpha<\gamma<\delta\) hence \(\alpha\leq s\) since \(\alpha\) is a limit point of \(p_\eta(\delta)\). Assuming \(\gamma=\xi_i\) and \(\delta= \eta_j\), 
we have \[\alpha\cap p_\xi(\xi_i) = \alpha\cap p_\eta(\eta_i) = \alpha\cap p_\eta(\eta_j),\] in which the first equality comes from the definition of \(f\) and the second comes from the coherence property of \(p_\eta\).
\end{proof}

\section{A nonspecial \(\omega_2\)-Aronszajn Tree}
We also need the following theorem from \Todorcevic~\cite{todorcevic1989special} (also see~\cite{todorcevic2007walks} Chapter~7). 

\begin{theorem}[\Todorcevic~\cite{todorcevic1989special}]\label{NonSpecialATree}
If \(\theta\geq \omega_2\) and there is a \(\Box(\theta)\)-sequence, then there is a nonspecial \(\theta\)-Aronszajn tree in the form of \(T(\rho_0^C)\) in which \(C\) is a nonspecial \(\Box(\theta)\)-sequence.
\end{theorem}

Shelah-Stanley~\cite{shelah1988weakly} has a similar result which states that \(\Box(\theta)\) implies the existence of a nonspecial \(\theta\)-Aronszajn tree. However, their construction is essentially different. The trees constructed in Shelah-Stanley~\cite{shelah1988weakly} admit \(\omega\)-ascent paths which make the nonspeciality of those trees robust. According to the historical remarks in~\cite{shelah1988weakly}, Laver first observed in 1970s that an \(\omega_2\)-Aronszajn tree with an ascent path cannot be special and later Baumgartner was able to construct such trees from \(\Box_{\omega_1}\) (see Devlin~\cite{devlin1983reduced}). However, as we will see soon, trees in the form of \(T(\rho_0^C)\) can be easily specialized.

\section{Specializing the \(\omega_2\)-Aronszajn Trees}
In this section, we show that the forcing axiom \(FA_{\omega_2}(\mathfrak{K})\) is strong enough to specialize all the trees in the form of \(T(\rho_0^C)\) for any nontrivial \(\Box(\omega_2)\)-sequence \(C\).

The following is a basic property of the sequence \(\rho_2\).
\begin{theorem}[\Todorcevic~\cite{todorcevic2007walks}]\label{rho2Unbounded}
  The following are equivalent for any C-sequence \(C = \langle C_\alpha : \alpha\in\Lim(\theta) \rangle\) on a regular uncountable cardinal \(\theta\) and the corresponding \(\rho_2\) sequence.
  \begin{enumerate}
    \item \(C\) is nontrivial.
    \item For every family \(A\) of \(\theta\) pairwisely disjoint finite subsets of \(\theta\) and every integer \(n\), there is a subfamily \(B\) of \(A\) of size \(\theta\) such that \(\rho_2(\alpha, \beta)>n\) for all \(\alpha\in a\), \(\beta\in b\) and \(a\neq b\) in \(B\).
  \end{enumerate}
\end{theorem}

The following lemma is essentially Lemma~6.3.3 of~\cite{todorcevic2007walks}. We include the proof here for completeness.

\begin{lemma}[\Todorcevic]\label{everySubsetContainsAnAntichain}
For any C-sequence \(C\) on \(\omega_2\), the tree \(T = T(\rho_2)\) has the property that every subset \(X\subset T\) of cardinality \(\omega_2\) contains an antichain of cardinality \(\omega_2\).
\end{lemma}
\begin{proof}
Consider a subset \(X\subset T\) of cardinality \(\omega_2\). We may assume that \(X\) is a level set and by replacing \(X = \{x_\alpha : \alpha\in\Gamma\}\) by a set lying inside its downwards closure, we may assume that the set consists of successor nodes of the tree. Let \(K\subset [\omega_2]^2\) be such that \[X = \{x_\alpha : \xi\in\Gamma\} = \{\rho_2(\cdot, \beta)\upharpoonright(\alpha+1) : \{\alpha,\beta\}\in K\}.\]
Shrinking \(X\) further, we may assume that pairs in \(K\) are pairwisely disjoint and \(\rho_2\) is constantly \(n\) on \(K\) since \(\rho_2\) only takes values of natural numbers. Apply Theorem~\ref{rho2Unbounded} we have \(K_0\subset K\) of cardinality \(\omega_2\) such that for all \(\{\alpha, \beta\}, \{\gamma, \delta\}\in K_0\) such that \(\alpha<\beta, \gamma<\delta\) and \(\alpha<\gamma\), we have \(\rho_2(\alpha, \gamma)>n\). Then \[X_0 = \{\rho_2(\cdot, \beta)\upharpoonright(\alpha+1) : \{\alpha,\beta\}\in K_0\}\] is an antichain in \(T\).
\end{proof}

We also need the following property of the tree \(T(\rho_2)\).

\begin{lemma}(CH)\label{rho2IsSufficientlyCoherent}
For any nontrivial \(\Box(\omega_2)\)-sequence \(\langle C_\alpha : \alpha\in\omega_2 \rangle\), the tree \(T = T(\rho_2)\) has the following property:

For any level sequence \(\langle p_\xi : \omega\rightarrow (T)_\xi \mid \xi\in\omega_2 \rangle\) such that every \(p_\xi\) is injective, there is a stationary \(\Gamma\subset \omega_2\) such that for any \(\xi<\eta\in\Gamma\) and any \(i, j\in\omega\), \[p_\xi(i)<p_\eta(i)\iff p_\xi(j)<p_\eta(j).\]
\end{lemma}
\begin{proof}
For any \(\delta\in\omega_2\) such that \(\cf \delta = \omega_1\) and any countable set \(X\subset\omega_2\) above \(\delta\), there exists a \(\delta'<\delta\) such that for any \(\alpha\in X\) and \(\xi\in tr(\delta, \alpha)\), either \(\delta\in \Lim C_\xi\) or \(C_\xi\cap\delta\subset \delta'\).

Clearly if \(\xi\in tr(\delta, \alpha)\) and \(\delta\in \Lim C_\xi\), the \(\xi\) must be the least element of \(tr(\delta, \alpha)\) above \(\delta\). Therefore, \(\delta'\) has the property that if \(\gamma\in[\delta', \delta)\), then the greatest element of \(tr(\gamma, \alpha)\) below \(\delta\) is \(\min(C_\delta\setminus\gamma)\) by the coherence of the square sequence.

Now if \(\alpha\in X\) and \(\gamma\in[\delta', \delta)\), and \(\xi = \min(tr(\delta, \alpha)\setminus(\delta+1))\), then
  \begin{equation*}
   \rho_2(\gamma, \alpha)=
    \begin{cases}
      \rho_2(\gamma, \delta) + \rho_2(\delta, \alpha) & \mbox{if \(\delta\in\Lim C_\xi\),}\\
      \rho_2(\gamma, \delta) + \rho_2(\delta, \alpha) -1 & \mbox{otherwise.}\\
    \end{cases}\end{equation*}

For each \(p_\xi\), let \(X_\xi = \{\alpha_\xi^i : p_\xi(i) = \rho_2(\cdot, \alpha_\xi^i)\upharpoonright\xi\}\). Apply the previous argument to every \(\xi\) and get an ordinal \(\xi'<\xi\). By shrinking the index set to a stationary \(\Gamma\subset\omega_2\), we may assume that for any \(\xi\in\Gamma\), \(\xi'<\xi\) is an constant ordinal \(s\) and \(p_\delta(i)\upharpoonright s\) is constantly \(t_i\). For any \(\xi\in\Gamma\) and \(\gamma\in[s, \xi)\), \[\rho_2(\gamma, \alpha_\xi^i) = \rho_2(\gamma, \xi) + \rho_2(\xi, \alpha_\xi^i) + \epsilon_\xi^i,\] in which \(\epsilon_\xi^i\) is either \(0\) or \(-1\).

By further shrinking \(\Gamma\), we may assume that neither \(\rho_2(\xi, \alpha_\xi^i)\) nor \(\epsilon_\xi^i\) depends on \(\xi\). Therefore, if \(\rho_2(\gamma,\alpha_\xi^i) = \rho_2(\gamma,\alpha_\eta^i)\) for some \(i\), it is also true for any \(i\), which completes the proof.
\end{proof}

\begin{lemma}(CH)
For any nontrivial \(\Box(\omega_2)\)-sequence \(\langle C_\alpha : \alpha\in\omega_2 \rangle\), there is a countably compact, well-met poset with \(\omega_2\)-Knaster poset which specializes the tree \(T(\rho_0)\)
\end{lemma}
\begin{proof}
Since there is a natural order-preserving mapping from \(T(\rho_0)\) to \(T(\rho_2)\), it suffices to specialize \(T = T(\rho_2)\).

Let \(P\) be the poset of all the conditions \(p\) satisfying
\begin{enumerate}
  \item \(p\) is a partial function from \(\omega_1\) to \(T^\omega\) such that \(\dom p\) is countable, and
  \item for every \(\alpha\in\dom p\), \(p(\alpha)\) is a countable antichain of \(T\),
\end{enumerate}
ordered by reversed inclusion.

It is clear that \(P\) is countably compact, well-met and it specializes the tree \(T\). It suffices to verify the \(\omega_2\)-Knaster property of \(P\).

Fix a sequence \(\langle p_\xi : \xi\in\omega_2 \rangle\). By CH and a standard \(\Delta\)-system argument, we may assume that the \(\dom p_\xi\)'s are constantly \(D\). For the sake of simplicity, we only verify the \(\omega_2\)-Knaster property for the case that \(D\) is a singleton \(\{d\}\) and identify \(p_\xi\) with \(p_\xi(d):\omega\rightarrow T\). In the general case that \(D\) is a countable subset of \(\omega_1\), the proof is the same since \(p_\xi(d)\) is already a countable set.

For \(p_\xi\) such that \(\cf \xi = \omega_1\), define \[s_\xi = \sup (\xi \cap (\{\hht(p_\xi(i)): i\in\omega\} \cup \{\Delta(p_\xi(i), p_\xi(j)) : i, j\in \omega\})).\] We have \(s_\xi < \xi\) and without loss of generality we may assume that \(s_\xi\) is constantly \(s\) for all \(\xi\). By futher shrinking the set, we may assume that \(p_\xi\) is a function \(p_\xi : \omega\rightarrow (T)_\xi\) and there is a sequence of nodes \(\langle t_i : i\in \omega\rangle\) from \((T)_{\hht(s)}\) such that for any \(\xi\), \(p_\xi(i)\upharpoonright \hht(s) = t_i\).

Now we can directly apply Lemma~\ref{rho2IsSufficientlyCoherent} and Lemma~\ref{everySubsetContainsAnAntichain} to this level set and this completes the proof.
\end{proof}

If the forcing axiom \(FA_{\omega_2}(\mathfrak{K})\) held, then for any \(\Box(\omega_2)\)-sequence \(C\), the tree \(T(\rho_0^C)\) must be special. On the other hand, since this axiom implies (S), hence \(\Box_{\omega_1}\), there exists such a tree which is nonspecial by Theorem~\ref{NonSpecialATree}, This is a contradiction and finishes the proof of Theorem~\ref{FAKnaster}.

\section{Conclusion}
As we already discussed, this inconsistency result is different from Shelah's Theorem~\ref{ccInconsistent}. It focuses on the chain condition rather than the well-metness condition. A related result can be found in \Lucke~\cite{lucke2017ascending}.

\begin{theorem}[\Lucke~\cite{lucke2017ascending}]\label{FAImpliesAWeaklyCompact}
  Let \(\kappa\) be an uncountable regular cardinal with \(\kappa = \kappa^{<\kappa}\). If \(\kappa^+\) is not weakly compact in \(L\), then there is a \(\kappa\)-closed, well-met partial order \(P\) satisfying the \(\kappa^+\)-chain condition such that FA\(_{\kappa^+}(P)\) fails.
\end{theorem}

On the other hand, Theorem~\ref{FAKnaster} in together with Theorem~\ref{ccInconsistent} sets a pretty tight bound on the generalization of Martin's Axiom to \(\omega_2\). If we compare this inconsistent axiom, FA\(_{\omega_2}(\mathfrak{K})\), with (S), theoretically speaking, there is an essential gap between the \(\omega_2\)-Knaster property and \(\omega_2\)-normality. For any \(\omega_2\)-sequence of forcing conditions, the pressing-down function requires the understanding of at least a club many of them, while Knaster property only requests an arbitary subset with cardinality \(\omega_2\). The different destinies of the corresponding forcing axioms also point to and result from that gap. However pragmatically, this gap is narrow in the sense that many proofs involving the Knaster property actually factorize through constructing pressing-down functions, which suggests that somehow Shelah's forcing axiom (S) is quite close to optimal.

\section{Acknowledgement}
The research of the first author is partially supported by grants from NSERC(455916) and CNRS(UMR7586).

The second author would like to thank Prof. Justin Moore for helpful conversation. The research of the second author is partially supported by NSF grant DMS-1854367.

\Addresses

\bibliographystyle{siam}
\bibliography{bibliography.bib}
\end{document}